\documentclass[a4paper,12pt,reqno]{amsart}
\usepackage{amssymb}
\usepackage[dvips]{graphicx}
\usepackage{hyperref}

\input epsf.tex
\newdimen\xsize
\newdimen\oldbaselineskip
\newdimen\oldlineskiplimit
\def\restorelineskip{\baselineskip=\oldbaselineskip%
\lineskiplimit=\oldlineskiplimit}
\def\putm[#1][#2]#3{
\hbox{\vbox to 0pt{\parindent=0pt%
\vskip#2\xsize\hbox to0pt{\hskip#1\xsize $#3$\hss}\vss}}}%
\long\def\Line#1{\hbox to \hsize{#1}}
\def\putt[#1][#2]#3{
\vbox to 0pt{\noindent\hskip#1\xsize\lower#2\xsize%
\vtop{\restorelineskip#3}\vss}}

\makeatletter
\def\xbig[#1]#2{{\hbox{$\m@th\left#2\vbox to#1\xsize{}%
\right.\n@space$}}}
\makeatother
\def\xlar[#1]#2{%
\smash{\mathop{ \hbox to #1\xsize{\leftarrowfill}}\limits^{#2}}}
\def\xrar[#1]#2{%
\smash{\mathop{ \hbox to #1\xsize{\rightarrowfill}}\limits^{#2}}}
\def\xline[#1]{\hbox to #1\xsize{\leaders\hrule\hfill}}

\DeclareFontFamily{U}{rsf}{\skewchar\font'177}%
\DeclareFontShape{U}{rsf}{m}{n}{<-6>rsfs5<6-8>rsfs7<8->rsfs10}{}%
\DeclareFontShape{U}{rsf}{b}{n}{<-6>rsfs5<6-8>rsfs7<8->rsfs10}{}%
\DeclareMathAlphabet\RSFS{U}{rsf}{m}{n}
\SetMathAlphabet\RSFS{bold}{U}{rsf}{b}{n}
\def\sf#1{{\mathsf{#1}}}

\def\slsf{\slshape \sffamily }

\def\msmall#1{\mathchoice{\hbox{\small$\displaystyle {#1}$}}{#1}{#1}{#1}}

\def\cc{{\mathbb C}}

\def\rr{{\mathbb R}}

\def\nn{{\mathbb N}}

\def\pp{{\mathbb P}}
  
\def\zz{{\mathbb Z}}
\def\ttt{{\mathbb T}}

\def\adyn{\sf{1}}

\def\span{\sf{span}}

\def\lim{\mathop{\sf{lim}}}
\def\limsup{\mathop{\sf{lim\,sup}}}

\def\min{\sf{min}}

\def\sup{\sf{sup}\,}

\def\v{{\mathrm{v}}}

\def\eps{\varepsilon}

\def\<{\langle}\let\la=\<
\def\>{\rangle}\let\ra=\>
  
\def\comp{\Subset}
\def\d{\partial}

\def\ddef{\mathrel{{=}\raise0.3pt\hbox{:}}}
\def\deff{\mathrel{\raise0.3pt\hbox{\rm:}{=}}}

\def\fraction#1/#2{\mathchoice{{\msmall{ #1\over#2}}}%
{{ #1\over #2 }}{{#1/#2}}{{#1/#2}}}
\def\norm#1{\left\Vert{#1}\right\Vert}

\def\le{\leqslant}

\def\longpoints{\leaders\hbox to 0.5em{\hss.\hss}\hfill \hskip0pt}
\def\stateskip{\smallskip}
\def\state#1. {\stateskip\noindent{\bf#1. }} 
\def\statep#1. {\stateskip\noindent{\bf#1 }} 
\def\proof{\state Proof. \2}

\def\Chi{\raise 2pt\hbox{$\chi$}}

\def\ie{\hskip1pt plus1pt{\sl i.e.\/,\ \hskip1pt plus1pt}}

\def\sli{{\sl i)} } 
\def\slii{{\sl i$\!$i)} }

\def\Chi{\raise 2pt\hbox{$\chi$}}

\let\phI=\phi\let\phi=\varphi\let\varphi=\phI
\let\cal=\mathcal

\def\calb{{\cal B}}
\def\calc{{\cal C}}
\def\cald{{\cal D}}
\def\cale{{\cal E}}
\def\calf{{\cal F}}
\def\calg{{\cal G}}

\def\calp{{\cal P}}

\def\cals{{\cal S}}

\def\calu{{\cal U}}
\def\calv{{\cal V}}

\def\calx{{\cal X}}
\def\caly{{\cal Y}}
\def\calz{{\cal Z}}

\def\eps{\varepsilon}

\def\comp{\Subset}
\def\d{\partial}

\def\1{{1\mkern-5mu{\rom l}}}

\def\ge{\geqslant}

\def\fraction#1/#2{\mathchoice{{\msmall{ #1\over#2}}}%
{{ #1\over #2 }}{{#1/#2}}{{#1/#2}}}

\def\le{\leqslant}

\newcommand{\2}{\thinspace}

\def\qed{\ \ \hfill\hbox to .1pt{}\hfill\hbox to .1pt{}\hfill $\square$\par}

\def\comment#1\endcomment{}

 \belowdisplayskip=\abovedisplayskip

\def\lineeqqno(#1){\hfill\llap{\vbox to 10pt%
{\vss\begin{align} \eqqno(#1)\end{align}\vss}}\vskip1pt}

\advance\headheight 1.2pt

\def\ShowwLLabel#1{}

\def\thechpt{\Roman{chpt}}

\def\newchapt[#1]#2{\newpage%
\refstepcounter{chpt}\setcounter{subsection}{0}%
\setcounter{thm}{0}\setcounter{defi}{0}%
\setcounter{rema}{0}\setcounter{exrc}{0}%
\renewcommand{\thesubsection}{\thechpt.\arabic{subsection}}%
\section*{\begin{center}\huge \bf Chapter \thechpt\\
#2 \end{center}}\label{#1}%
\ \smallskip%
\markboth{Chapter \thechpt}{#2}%
}
\def\newsect[#1]#2{\refstepcounter{section}\setcounter{equation}{0}%
\renewcommand{\thesubsection}{\arabic{section}.\arabic{subsection}}%
\section*{\arabic{section}.
#2}\vspace{-20pt}\label{#1}\vspace{20pt}%
\markboth{Section \arabic{section}}{#2}}

\def\newlect[#1]#2{\refstepcounter{section}%
\renewcommand{\thesubsection}{\arabic{section}.\arabic{subsection}}%
\section*{Lecture \arabic{section}\\
#2}\label{#1}%
\markboth{Lecture \arabic{section}}{#2}}

\def\newprg[#1]#2{\refstepcounter{subsection}%
\subsection*{{\thesubsection.\ #2}} \label{#1}%
}

\def\newappx[#1]#2{%
\refstepcounter{appx}\setcounter{section}{0}%
\renewcommand{\thesubsection}{A\arabic{appx}.\arabic{subsection}}%
\section*{Appendix \arabic{appx}\\ #2}
\label{#1}%
\markboth{Appendix A\arabic{appx}}{#2}
}

\newtheorem{thm}{Theorem}[section]
   \def\newthm#1{\begin{thm}\label{#1}}

\newtheorem{nnthm}{Theorem}
   \def\newthm#1{\begin{nnthm}\label{#1}}

\newtheorem{lem}{Lemma}[section]
   \def\newlemma#1{\begin{lem} \label{#1}}

\newtheorem{prop}{Proposition}[section]
   \def\newprop#1{\begin{prop}\label{#1}}

\newtheorem{nnprop}{Proposition}
   \def\newprop#1{\begin{nnprop}\label{#1}}

\newtheorem{corol}{Corollary}[section]
   \def\newcorol#1{\begin{corol} \label{#1}}

\newtheorem{nncorol}{Corollary}
   \def\newcorol#1{\begin{nncorol} \label{#1}}

\newtheorem{defi}{Definition}[section]
   \def\newdefi#1{\begin{defi} \label{#1}\rm }

\newtheorem{nndefi}{Definition}
   \def\newdefi#1{\begin{nndefi} \label{#1}\rm }

\newtheorem{exmp}{Example}[section]
   \def\newexmp#1{\begin{exmp} \label{#1}\rm }

\newtheorem{exrc}{Exercise}
   \def\newexrc#1{\begin{exrc} \label{#1}\rm }

\newtheorem{rema}{Remark}[section]
   \def\newrema#1{\begin{rema} \label{#1}\rm }

\newtheorem{nnrema}{Remark}
   \def\newrema#1{\begin{rema} \label{#1}\rm }

\def\eqqno(#1){\label{(#1)}}
\def\eqqref(#1){(\ref{(#1)})}

\usepackage[active]{srcltx}   

\def\el2{\sf{L^2}}
\def\el{\sf{l}}

\title{On Hilbert-Hartogs manifolds}
\author{M. Anakkar*, A. Zagorodnyuk**}
\address{*Univ.Lille,CNRS, UMR 8524-Laboratoire Paul Painlevé, F-59000 Lille, France}
\email{mohammed.anakkar@univ-lille.fr}
\address{**
Precarpathian National University, Ivano-Frankivsk, Shevchenka
str. 57, 76018, Ukraine.} 
\email{andriyzag@yahoo.com}

\thanks{* Supported by the Labex Cempi ANR-11-LABX-0007-01}
\thanks{** Partially supported by French-Ukrainian CNRS/NANU, Grant no.: 24013.
*** Supported by Grant DFFD of Ukraine F35/531-2011}
\date{\today}

\begin{document}

\begin{abstract}
We  prove some  extension results for holomorphic mappings with values in complex Hilbert manifolds.
\end{abstract}
\maketitle

\setcounter{tocdepth}{1}
\tableofcontents

\newsect[INT]{Introduction}

Throughout this paper by capital latin letters, like $S,X,Y$, we shall denote finite dimensional complex 
manifolds, by calligraphic letters, like $\cals, \calx , \caly$ - complex Hilbert manifolds. Our Hilbert manifolds 
are modeled over $\el^2$ and are supposed to be second countable. Recall that a $q$-concave Hartogs figure in $\cc^{q+1}$
is the following domain 
\begin{equation}
\eqqno(hart-11)
H_q^1(r) = (\Delta^q \times \Delta(r) ) \cup (A_{1-r,1}^q\times \Delta).
\end{equation}
Here $\Delta (r)=\Delta_r$ denotes the disk of radius $r$ in $\cc$ centered at zero, $\Delta$ the unit disk,
$\Delta^q(r)$ the polydisk in $\cc^q$, $A_{1-r,1}^q = \Delta^q\setminus \bar \Delta_{1-r}^q$ the ring domain. 
$r>0$ is a small positive number, its precise value is usually irrelevant: if something holds for some $0<r<1$ 
then the same usually holds for all other $0<r'<1$. We say that a comlex manifolds $X$ is $q$-Hartogs if every 
holomorphic mapping $f:H_q^1(r) \to X$ extends to a holomorphic mapping $\tilde f:\Delta^{q+1}\to X$ of the unit 
$(q+1)$-disk to $X$. If the same holds for a complex Hilbert manifold $\calx$ we call $\calx$ $q$-Hilbert-Hartogs. 
Let $D$ be a domain in a complex  manifold $X$ with $\mathcal{C}^2$ boundary. $D$ is called $q$-pseudoconcave at a
boundary point $p$ if the Levi form of a defining function of $D$ near $p$ has at least $q$ negative eigenvalues on $T^c_p\d D$. 
A complex Hilbert manifold $\calx$ satisfies the $q$-Levi extension condition if for any domain $D \subset \mathbb{C}^n$
with $\mathcal{C}^2$ boundary and such that $\partial D$  is q-pseudoconcave at $p$, any holomorphic map $f: D \rightarrow \calx$ 
extends holomorphically to a neighborhood of p. Our main result in this paper is the following

\begin{nnthm}
\label{hart-levi}
A complex Hilbert manifold $\calx$ satisfies the $q$-Levi extension condition if and only if 
$\calx$ is $q$-Hartogs.
\end{nnthm}

The proof of this and other extension results for complex Hilbert manifolds relies on one 
technical statement which is a week version of Royden's lemma for Hilbert manifolds, and we are going to 
formulated it now. Recall that an open subset $\calv\subset \calx$ is 
called $1$-complete if it admits a strictly plurisubharmonic exhaustion function see Definition \ref{psh-f} below. 
By an analytic disk we mean a holomorphic mapping $\phi$ of some neighborhood of $\bar\Delta$ into $\calx$. By an 
analytic $q$-disk in $\calx$ we mean a holomorphic mapping of a neighborhood of the closure $\bar D$ of a bounded 
strictly pseudoconvex domain $D\subset \cc^q$ to $\calx$.    The technical statement we ment is the following.

\begin{nnthm}
\label{royden-w}
Let $\phi : \bar D \to \calx $ be an imbedded analytic $q$-disk in a complex Hilbert manifold $\calx$. 
Then $\phi (\bar D) $ has a fundamental system of $1$-complete neighborhoods.
\end{nnthm}

\begin{nnrema} \rm
\label{polydisk}
This theorem is valid also for such $D$ as polydisc $\Delta^q$ or, a product $B^{q-k}\times \Delta^k$ of a ball 
with polydisk. The latter two are not smoothly bounded, but since $\phi$
is supposed to be defined in a neighborhood of $\bar D$, we can find a smoothly bounded s.p.c. $D_1\supset \bar D$
contained in this neighborhood and then apply our theorem to $D_1$.
\end{nnrema}

\smallskip\noindent{\slsf The structure of the paper.} 

\smallskip\noindent{\slsf 1.} In section \ref{ROYD}, after  some preliminaries, we prove the existence of 
$\adyn$-complete neighborhoods, \ie Theorem \ref{royden-w}.

\smallskip\noindent{\slsf 2.} Section \ref{HART} is devoted to the general study of Hilbert-Hartogs manifolds.
In particular we prove Theorem \ref{hart-levi} there as well as some other results.

\smallskip\noindent{\slsf Asknowledgments.} Authors would like to express their thanks to S. Ivashkovich
and L. Blanc-Centi for the helpfull discussions on the subject of this paper.

\newsect[ROYD]{Existence of complete neighborhoods}

\newprg[ROYD.hilb]{Hilbert manifolds}

The notation $\el^2$ throughout this paper stands for the Hilbert
space of sequences of complex numbers $z=\{z_k\}_{k=1}^{\infty}$
such that $||z||^2:=\sum_k|z_k|^2<\infty $ with the standard
Hermitian scalar product $(z,w) = \sum_kz_k\bar w_k$ and standard
basis $\{ e_1,e_2,...\}$.  The
$q$-dimensional complex linear space $\cc^q$ will be often
identified with the subspace $\el^2_q = \span\{e_1,...,e_q\} \subset
\el^2$ and we shall often write (with some ambiguity) $\el^2 =
\cc^q\oplus \el$, where $\el$ is the orthogonal complement to $\cc^q$ in $\el^2$
and use coordinates $z'=(z_1,...,z_q)$ for $\cc^q$ and
$z^{''}=\{z_{q+1},z_{q+2},...\}$ for $\el$. Moreover, we shall think
about $\el$ as of a copy of $\el^2$ itself and henceforth using
notations like $B^q\times B^{\infty}$ to denote an obvious subset of
$\el^2$.

\smallskip Hilbert manifolds of this paper are modeled over $\el^2$, \ie they are
Hausdorff topological spaces locally homeomorphic to open subsets of $\el^2$ with biholomorphic 
transition mappings. In addition they are supposed to be second countable.
The existence of a partition of unity on such Hilbert manifolds subordinated to a given open 
covering can be proved following \cite{Mu} as follows.

\smallskip Let $\calx$ be a Hilbert manifold and $(\Omega_\alpha,\phi_\alpha)_{\alpha \in A}$ be
an atlas. Since $\calx$ is second countable, one can suppose that $A$ is countable. 
Consider the subspace $\mathbb{Q}^{\mathbb{N}}_0$ of $\el^2$ which consists of finite sequences 
of rational numbers. For every $x \in \mathbb{Q}^{\mathbb{N}}_0 \cap 
\phi_\alpha(\Omega_\alpha)$ there exists $\epsilon_x$ such that $B(x,2\epsilon_x) \subset 
\phi_\alpha(\Omega_\alpha)$. Then $\{\phi_\alpha^{-1}(B(x,\epsilon_x)) \ | \ x \in  
\mathbb{Q}^{\mathbb{N}}_0 \cap \phi_\alpha(\Omega_\alpha)\}$ is a countable covering of 
$\Omega_\alpha$ and, since $A$ is countable, one obtains a countable 
covering $\{B_n\}_{n\in \nn}$ of $\calx$ subordinated to $\{\Omega_\alpha \}_{\alpha \in A} $. 
By the axiom of choice there is a function 
$\tau: \mathbb{N} \to A$ such that $B_n \subset \Omega_{\tau(n)}$.
 Construct a sequence of smooth functions $0 \leqslant f_n \leqslant 1$ such that
\begin{equation}
f_n(x) = 
     \begin{cases}
       1 & \quad{\text{if}} \ ||\phi_{\tau(n)}(x)-x_n||<\epsilon_{x_n}\\
       0& \quad{\text{if}} \ ||\phi_{\tau(n)}(x)-x_n||>2\epsilon_{x_n}.\\ 
     \end{cases}
\end{equation} 
Then define another sequence $(g_n)_n$ by
\begin{equation}
\begin{cases}
       g_1 = & f_1, \\
       g_n = & f_n \ \prod_{j=1}^{n-1}(1-f_j) \ \quad{\text{if}} \ n \geqslant 2. \\
\end{cases}
     \end{equation} 
 Clearly $0 \leqslant g_n \leqslant 1$ on $\calx$ and $supp(g_n) \subset \Omega_{\tau(n)}$ for every $n$.
 Furthermore, one can prove by induction that 
 \begin{equation}
 \eqqno(sum-gn)
 g_1 + ... + g_n = 1 - \prod_{j=1}^{n}(1-f_j).
 \end{equation}
 Remark that

\begin{itemize}
 \item  From \eqqref(sum-gn) it follows that  $g_1 + ... + g_n =1$ on $B_n$ 
 because $f_n=1$ on $B_n$. 

\smallskip
 \item From the definition of $(g_k)$, we see that $g_k=0$ on $B_n$ for every $k>n$.
\end{itemize}

The first item guarantees that $\sum_{n \in \mathbb{N}}{g_n(x)}=1$ for every $x\in  \calx$ and the 
second one gives that $(g_n)$ is locally finite on $\calx$ . Construct then the partition of unity 
$\{\psi_\alpha\}_\alpha$  subordinated to $\{\Omega_\alpha\}_{\alpha \in A}$ by:
 \begin{equation}
 \psi_\alpha=
 \begin{cases}
       \sum_{\tau(n)=\alpha}{g_n} & \quad{\text{if}} \ \alpha \in \tau(\mathbb{N}),\\
       0 & \quad{\text{otherwise.}} \\ 
     \end{cases} 
\end{equation}

\begin{rema} \rm
From the partition of unity subordinated to the atlas of the Hilbert complex maniflod $\calx$ one can 
construct a Hermitian metric $g$ on $\calx$ by patching the Hermitian metrics on charts.
\end{rema}

\smallskip $\el^2$-valued holomorphic mapping on an open subset $\cald$ of $\el^2$ is a mapping
$f:\cald\to \el^2$ which is Fr\'echet differentiable at all points of $\cald$
or, equivalently, which in an a neighborhood of every point $z^0\in
\cald$ can be represented by a convergent power series
\begin{equation}
\eqqno(hol-map)
f(z) = f(z^0) + \sum\limits_{n=1}^{\infty}P_n(z-z^0),
\end{equation}
where $P_n$ are {\slsf continuous homogeneous} polynomials of degree $n$ satisfying
\begin{equation}
\eqqno(cauchy)
\limsup\limits_{n\to\infty}\norm{P_n}^{\frac{1}{n}} = \frac{1}{r} <\infty .
\end{equation}
The radius of convergence $r(z^0)$ given by \eqqref(cauchy) can be less than
the distance $d(z^0,\d \cald)$ from $z^0$ to the boundary of $\cald$. But if $f$
is bounded on $\cald$ then $r(z^0)\ge d(z^0,\d \cald)$, see \cite{Mu} for more details.
The norm of a continuous $n$-homogeneous polynomial is defined as
\begin{equation}
\eqqno(norm-n)
\norm{P} \deff \sup\{\norm{P(x)}:\norm{x}\le 1\}.
\end{equation}
Cauchy-Hadamard formula \eqqref(cauchy) guarantees the uniform and
absolute (\ie normal) convergence of power series \eqqref(hol-map) on every ball
$B(z^0,r')$ with $r'<r(z^0)$. It is important to note that the space
$\calp_n$ of continuous, homogeneous polynomial mappings of degree
$n$ from $\el^2$ is a Banach space with respect to the norm
\eqqref(norm-n).

\medskip We shall need a version of Grauert's theorem due to Bungart, see \cite{Bu}
Theorem 3.2 and/or \cite{Le} Theorem 9.2, which, combined with the Theorem of Kuiper about 
contractibility of the group of invertible operators in infinite dimensional Hilbert space, 
see \cite{Ku}, gives the following

\begin{thm}
\label{grauert} Let $D$ be a Stein manifold and $\cale$ a
holomorphic Hilbert vector bundle over $D$. Then $\cale$ is
holomorphically trivial and, moreover,  if $\calu = \{U_{\alpha}\}$ is a
locally finite Stein covering of $D$ then for every cocyle $f\in
Z^1(\calu ,\calg )$ there exists a cochain $c\in C^0(\calu
,\calg )$ such that $\delta (c) = f$.
\end{thm}
Here a Hilbert bundle is understood as being locally isomorphic to
products $U\times \el^2$ and $\calg $ is a sheaf of holomorphic mappings with values in
the group of invertible operators on $\el^2$.

\newprg[ROYD.loc]{Existence of $\adyn$-Complete Neighborhoods}
We are going to prove Theorem \ref{royden-w} from the introduction. This will be done 
following the original idea of H. Royden from \cite{Ro}. Let $\calv$ an open set of $\el^2$
and  $u \in \mathcal{C}^2(\calv, \rr)$.

\begin{defi}
\label{psh-f}
We say that $u$ is strictly plurisubharmonic if
\[
\sum_{k,l=1}^{\infty}{\frac{\partial^2u}{\partial z_k \partial\bar{z}_l}(z_0)v_k\bar{v}_l} 
\geqslant c(z_0)||v||^2 \quad \forall v\in \el^2,
\]
with $c(z_0)$ positive and continuous on $\calv$.
\end{defi}

We say that (continuous) $u:\calv\to [0, t_0)$ is an exhaustion function of $\calv$ if 

\begin{itemize}
 \item for every $t<t_0$ one has that $\overline{u^{-1}\left([0,t)\right)}\subset \calv$,
 
 \item for all $0\le t_1<t_2<t_0$ one has that $\overline{u^{-1}\left([0,t_1)\right)}\subset
 u^{-1}\left([0,t_2)\right)$,
 
 \item $\calv = u^{-1}\left([0,t_0)\right)$.
\end{itemize}

\begin{defi}
\label{1-compl}
$\calv$ is said to be $1$-complete if it admits a strictly plurisubharmonic exhaustion function.
\end{defi}

Recall that by an analytic $q$-disk in a complex Hilbert manifold $\calx$ we
understand a holomorphic mapping $\phi $ of a neighborhood of a
closure of a relatively compact strongly pseudoconvex domain $D\comp
\cc^q$ into $\calx$. The image $\phi (\bar D)$ we shall denote by
$\Phi$. Therefore, by saying that a $q$-disk $\Phi =\phi (\bar D)$
is analytic we mean that $\phi $ holomorphically extends to a
neighborhood of $\bar D$. Saying that $\Phi$ is imbedded we mean
that $\phi$ is an imbedding in a neighborhood of $\bar D$.

\smallskip  We choose $r>0$ such that $\phi $ extends as an imbedding to a 
$2r$-neighborhood of $D$, \ie to $D^{2r} \deff \{z\in \cc^q: d(z,D)< 2r\}$.

\begin{lem}
\label{redress}
Let $\phi : B^q(a,\eps)\to \calx$ be a holomorphic map of a ball centered
at $a\in\cc^q$ into a complex Hilbert manifold $\calx$ 
such that $d_a\phi:\cc^q\to T_{b}\calx$ is injective. Then one can find a
coordinate chart $(V, h)$ in a neighborhood of $b=\phi (a)$ such that:

\sli $V$ is mapped by $h$ onto a neighborhood $V'$ of the point $(a, 0)\in \el^2$
with $h(b) = (a, 0)$. If $z = ( z_1,..., z_q),  w =
( w_{q+1},,,,)$ are standard coordinates in $\el^2$ then $V' = \{( z,
 w)\in\el^2: || z-a||<\delta , || w||<\delta\}$ for an
appropriate $\delta >0$.

\slii Mapping $h\circ \phi $ from $U :=\phi^{-1}(V)$ to $V'$ in these coordinates is
given by $(h\circ \phi )(z_1,...,z_q)=(z_1,...,z_q,0,...)$.
\end{lem}
\proof Take some coordinate chart $(V_1,h_1)$ in a neighborhood of
$b$ in $\calx$. Choose a frame $\{ \v_i\}_{i=1}^{\infty}$ in $\el^2$ in
such a way that:

\smallskip{\slsf (a)} $d_{a}\phi (e_i) = \v_i, i=1,...q$, for the
standard basis $e_1,...,e_q$ of $\cc^q$;

\smallskip{\slsf (b)} $\span \left\lbrace \v_1,...,\v_q\right\rbrace$
is orthogonal to $\span \left\lbrace \v_{q+1},...\right\rbrace$.

\medskip Let $z' = (z'_1,...,z'_q), w' = (w_{q+1}',...)$ be affine coordinates
in $\el^2$ which correspond to the frame $\{\v_i\}$ and such that
$h_1(b) = (a, 0)$ in these coordinates. Write $h_1\circ \phi = (\phi_1,
\phi_2)$ in coordinates $(z',w')$. One has
$\frac{\d(z'_1,...,z'_q)}{\d(z_1,...,z_q)}(0) = \adyn_q$ due to (a).
By implicit function theorem there exist neighborhoods $U\ni a$ in
$\cc^q$ and $V_2\ni a$ in $L = {\span}\{\v_1,...,\v_q\}$ such
that $\phi_1:U \to V_2$ is a biholomorphism. $U$ can be taken of the
form $\{||z - a|| < \delta \}$ for an appropriate $\delta  >0$.

\smallskip Therefore $\phi (U)$ is a graph $w' = \psi
(z')$ over $V_2$. Make a coordinate change $h_2: (z^{'} ,w^{'}) \to (z^{''}=z^{'},
w^{''} = w^{'} - \psi(z^{'}))$ to get that $(h_2\circ h_1)\circ \phi$ has the
form $z \to (\phi_1(z), 0)$. Finally make one more coordinate
change $(z^{''},w^{''})\to ( z = \phi_1^{-1}(z^{''}),
 w = w^{''})$ to get the final chart $(V,h)$ with
$(h\circ \phi )(z) = (z,0)$. $V'$ can be taken to
have the form $U\times V_3$ where $V_3=\{w,
|| w||<\delta \}$ for an appropriate $\delta >0$ (for that one might need
to shrink $V_2$ and therefore $U$). Lemma \ref{redress} is proved.

\smallskip\qed

\begin{rema} \rm
$U$ was chosen to be a ball. Note that it can be chosen to be a cube as well.
Remark also that $V'$ can be taken in the form of a product $V' = U\times V{''}$,
where $V{''}$ is a Hilbert ball ($V^{''} = V_3$ in the notations of the proof
of Lemma \ref{redress}).
\end{rema}

\newprg[ROYD.triv]{Trivializaion of the infinitesimal neighborhood} Cover $\phi (\overline{D^r})$ 
with a finite collection of
coordinated neighborhoods $\{(V_{\alpha}, h_{\alpha})\}_{\alpha
=1}^N$ with centers at $a_{\alpha}$ as in Lemma \ref{redress}. Denote by
$z_{\alpha},w_{\alpha}$ the corresponding coordinates in
$V'_{\alpha} = U_{\alpha}\times V^{''}_{\alpha}\subset \el^2$.
Note that on this stage $z_{\alpha}$
glues to a global coordinate $z$ on $D$. Denote by $J_{\alpha ,
\beta}$ the Jacobian matrix of the coordinate change $(z_{\alpha},
w_{\alpha}) = (h_{\alpha}\circ h_{\beta}^{-1})(z_{\beta},
w_{\beta})$. Since $h_{\alpha}\circ h_{\beta}^{-1}(z, 0) = (z, 0)$
we see that

\begin{equation}
\eqqno(jak0)
J_{\alpha ,\beta} (z, 0) =
\begin{pmatrix}
I_q & A_{\alpha ,\beta}(z)\cr
0   & B_{\alpha ,\beta}(z)
\end{pmatrix}.
\end{equation}
By construction the operator valued functions $B_{\alpha ,\beta}$
form a multiplicative cocycle, \ie they are the transition functions
of an appropriated vector $\el^2$-bundle over $\overline{D^r}$ (the
normal bundle to $\phi (\overline{D^r})$). Indeed
\begin{equation}
\eqqno(jak1)
\begin{pmatrix}
I_q & A_{\alpha ,\beta}(z)\cr
0   & B_{\alpha ,\beta}(z)
\end{pmatrix}
\cdot
\begin{pmatrix}
I_q & A_{\beta ,\gamma}(z)\cr
0   & B_{\beta ,\gamma}(z)
\end{pmatrix} =
\begin{pmatrix}
I_q & A_{\beta ,\gamma}(z)+ A_{\alpha ,\beta}(z)B_{\beta ,\gamma}(z)\cr
0   & B_{\alpha ,\beta}(z)B_{\beta ,\gamma}(z)
\end{pmatrix} =
\begin{pmatrix}
I_q & A_{\alpha ,\gamma}(z)\cr
0   & B_{\alpha ,\gamma}(z)
\end{pmatrix}
\end{equation}
This bundle is trivial by Theorem \ref{grauert} and therefore one can find 
holomorphic operator
valued functions $B_{\alpha}: U_{\alpha} \to End (\el^2)$ such that
$B_{\alpha , \beta} = B_{\alpha}\circ B_{\beta}^{-1}$ on $U_{\alpha}
\cap U_{\beta}$. Make a coordinate change in $V_{\alpha}$ as follows:
$\tilde z_{\alpha} = z_{\alpha}$ (\ie still $\tilde z_{\alpha} = z$
for all $\alpha$)  and $\tilde w_{\alpha} = B_{\alpha}
(z_{\alpha})^{-1}w_{\alpha}$. Then in new coordinates the Jacobian
matrix of the coordinate change (when restricted to $w_{\beta} = 0$)
will have the form

\begin{equation}
\eqqno(jak2)
J_{\alpha ,\beta} (z; 0) =
\begin{pmatrix}
I_q & A_{\alpha ,\beta}(z)\cr
0   & I_{\infty}
\end{pmatrix} ,
\end{equation}
with some other $A_{\alpha ,\beta}(z)$. Notations for coordinates $(z_{\alpha}, w_{\alpha})$ and charts
$(V'_{\alpha}=U_{\alpha}\times V^{''}_{\alpha},
h_{\alpha})$ will be not changed on this stage. 


\begin{figure}[h]
\centering
\includegraphics[width=4.0in]{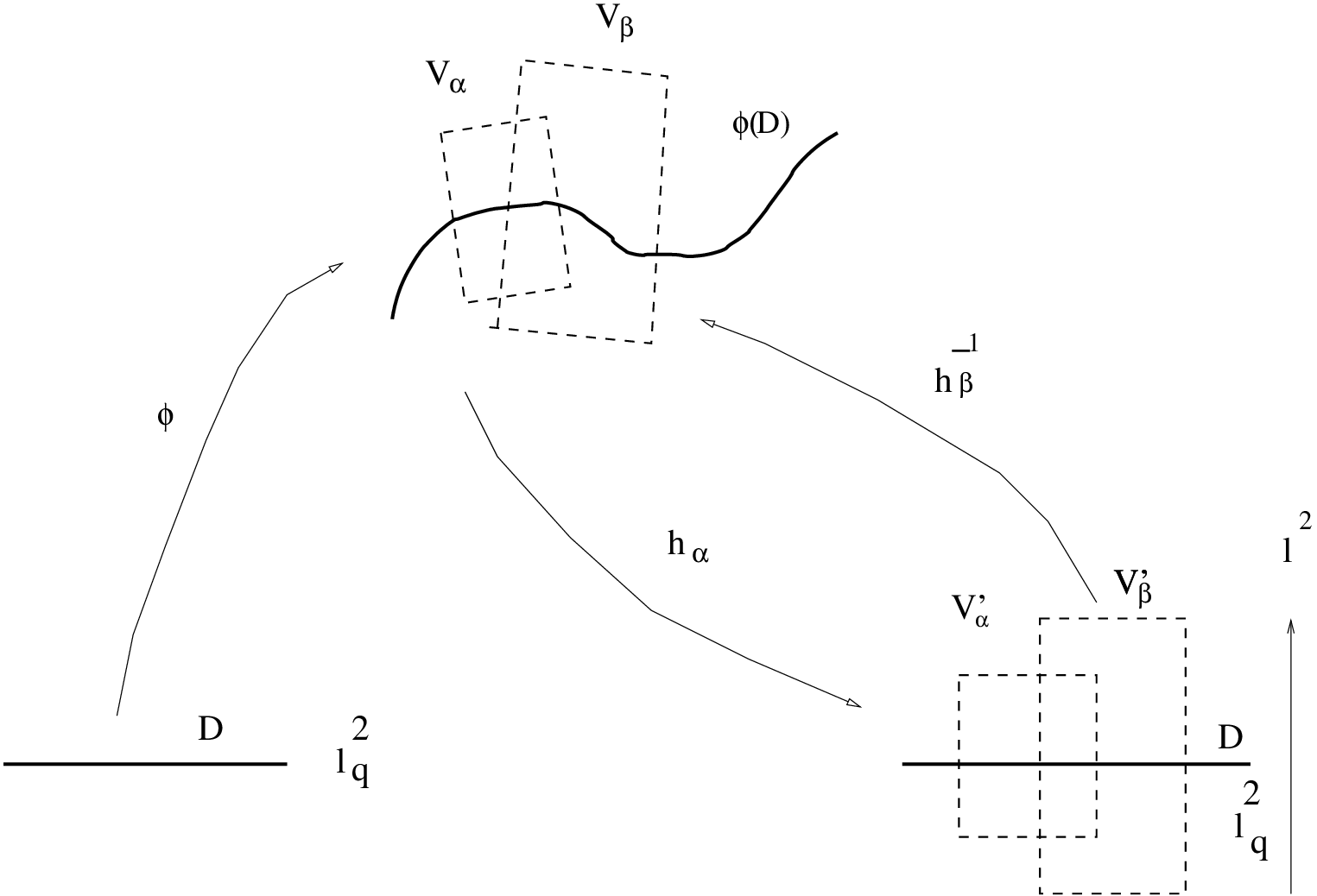}
\caption{$z=(z_1,...,z_q)$ serves as common coordinate for all charts on this step.}
\label{pic1}
\end{figure}

\newprg[ROYD.tr-hi]{Trivializaion to higher orders.}
Transition mappings between new coordinate charts  have the form

\begin{equation}
\eqqno(chng1)
\begin{cases}
z_{\alpha} = z_{\beta} + \sum\limits_{n = 1}^{\infty} A^n_{\alpha ,
\beta }(w_{\beta}),\cr
w_{\alpha} = w_{\beta} + \sum\limits_{n = 2}^{\infty} B^n_{\alpha , \beta }
(w_{\beta}),
\end{cases}
\end{equation}
where $A^n_{\alpha , \beta}$ (resp.  $B^n_{\alpha , \beta }$) are
holomorphic functions of $z_{\beta}=z\in U_{\alpha}\cap U_{\beta}$ with values in
the spaces of continuous $n$-homogeneous vector valued polynomials on $w_{\beta}$.
We shall use the standard notations, see \cite{Mu}, as follows.
For an $n$-homogeneous polynomial $B^n$ we denote as $\hat B^n$ the (unique) symmetric
$n$-linear form defining $B^n$, \ie such that $B^n(w) = \hat B^n(\underbrace{w,...,w}_{n-times}) =: \hat Bw^n$.
According to the same standard notations when writing $\hat B^nu^l v^k$ one means
$\hat B^n(\underbrace{u,...,u}_{l-times},\underbrace{v,...,v}_{k-times})$, $n=l+k$. Our goal in
what follows is to ged rid of terms $B^n_{\alpha , \beta}$. Remark that on $U_{\alpha}\cap U_{\beta}
\cap U_{\gamma}$ one has
\[
w_{\alpha} = w_{\beta} + \sum\limits_{n = 2}^{\infty} B^n_{\alpha , \beta }(w_{\beta}) =
w_{\gamma} +  \sum\limits_{n = 2}^{\infty} B^n_{\beta , \gamma }(w_{\gamma}) +
\]
\[
 + \sum\limits_{n = 2}^{\infty} B^n_{\alpha , \beta }\left(w_{\gamma} +
\sum\limits_{m = 2}^{\infty} B^m_{\beta , \gamma }(w_{\gamma})\right) =
w_{\gamma} + \sum\limits_{n = 2}^{\infty} B^n_{\alpha , \gamma }(w_{\gamma}).
\]
Therefore
\begin{equation}
\eqqno(rel-b)
\sum\limits_{n = 2}^{\infty} B^n_{\alpha , \gamma }(w_{\gamma}) =
 \sum\limits_{n = 2}^{\infty} B^n_{\beta , \gamma }(w_{\gamma})  + \sum\limits_{n = 2}^{\infty}
B^n_{\alpha , \beta }\left(w_{\gamma} + \sum\limits_{m = 2}^{\infty} B^m_{\beta , \gamma }(w_{\gamma})\right).
\end{equation}
This gives for every degree $N\ge 2$ of homogenuity the following {\slsf finite} relation
between homogeneous polynomials of degree $N$
\begin{equation}
\eqqno(rel-bN)
B^N_{\alpha , \gamma}(w_{\gamma}) = B^N_{\beta , \gamma}(w_{\gamma}) + \sum\limits_{(m-1)k+n=N}
C_{n-k}^k\hat B^n_{\alpha , \beta}w_{\gamma}^{n-k} B^m_{\beta , \gamma}(w_{\gamma})^k.
\end{equation}
In the right hand side
of \eqqref(rel-bN) we have the term $B^N_{\beta , \gamma}(w_{\gamma})$ and terms in the sum with $k\ge 0$.
Each of this terms should have the wright degree of homogenuity equal to $N$, which gives only a finite number of possibilities
for $k,m$ and $n$. Remark that if $N=2$ then, since $n, m\ge 2$, the only possibility in \eqqref(rel-bN)
is $k=0$ and therefore we get a cocyle condition
\begin{equation}
\eqqno(coc-b2)
B^2_{\alpha , \gamma}(w_{\gamma}) = B^2_{\beta , \gamma}(w_{\gamma}) +  B^2_{\alpha , \beta}(w_{\gamma}).
\end{equation}

\begin{rema} \rm
It is worth to point out at this stage that we do not have the cocycle condition for higher $N$-s, for them
we have only relation \eqqref(rel-bN).
\end{rema}

Now solve the additive Cousin problem  $B^2_{\alpha , \beta} = B^2_{\beta} - B^2_{\alpha}$ for the acyclic
covering $\{U_{\alpha} \}$ of $D$ and make the following (quadratic in $w$) change of variables
\begin{equation}
\eqqno(chng-b2)
\begin{cases}
\tilde z_{\alpha} = z_{\alpha}\cr
\tilde w_{\alpha} = w_{\alpha} + B^2_{\alpha}(w_{\alpha}).
\end{cases}
\end{equation}

It is not difficult to check that relation \eqqref(chng1) becomes now
\begin{equation}
\eqqno(chng2)
\begin{cases}
\tilde z_{\alpha} = \tilde z_{\beta} + \sum\limits_{n = 1}^{\infty} \tilde A^n_{\alpha , \beta } (\tilde w_{\beta}),\cr
\tilde w_{\alpha} = \tilde w_{\beta} + \sum\limits_{n = 3}^{\infty} \tilde B^n_{\alpha , \beta }(\tilde w_{\beta})
\end{cases}
\end{equation}
with an appropriate $\tilde A^n_{\alpha , \beta}$ and $\tilde B^n_{\alpha , \beta}$.

\medskip Now we can trivialize the normal bundle in all orders. 
From \eqqref(rel-bN) we see that if $B^n_{\alpha , \beta} =0$ for all $n<N$ and all
$\alpha, \beta$ then $\{B^N_{\alpha, \beta}\}$ is an additive cocyle, which can be resolved as
$B^N_{\alpha , \beta} = B^N_{\beta} - B^N_{\alpha}$ and then the coordinate change

\begin{equation}
\eqqno(chng-bN)
\begin{cases}
z_{\alpha} = z_{\beta}\cr
\tilde w_{\alpha} = w_{\alpha} + B^N_{\alpha}(z)(w_{\alpha}).
\end{cases}
\end{equation}
will give us coordinates in which the transition functions will take the form
\begin{equation}
\eqqno(chngN)
\begin{cases}
z_{\alpha} = z_{\beta} + \sum\limits_{n = 1}^{\infty} A^n_{\alpha ,
\beta } (w_{\beta}),\cr
w_{\alpha} = w_{\beta} + \sum\limits_{n = N+1}^{\infty} B^n_{\alpha , \beta }(w_{\beta}).
\end{cases}
\end{equation}
$\cc^q$-valued functions $A^n_{\alpha , \beta }(z_{\beta})$ are
holomorphic in $z\in U_{\alpha}\cap U_{\beta}$
with values in the space of continuous homogeneous polynomials of degree $n$ in $w$-s. They
satisfy the additive cocycle condition (for every fixed $n\le N$).

\smallskip For every fixed $n\le N$ we solve $A^n_{\alpha ,
\beta}=A^n_{\beta} - A^n_{\alpha}$  and in every chart $V'_{\alpha}$ make the coordinate
change:

\begin{equation}
\eqqno(chng6)
\begin{cases}
\tilde z_{\alpha} = z_{\alpha} - \sum\limits_{n=1}^{N}A^n_{\alpha}(z_{\alpha})(w_{\alpha})\cr
\tilde w_{\alpha} = w_{\alpha}.
\end{cases}
\end{equation}
We get that in new coordinates the coordinate changes vanish to order $N$, \ie have the form
\begin{equation}
\eqqno(chngN+1)
\begin{cases}
z_{\alpha} = z_{\beta} + \sum\limits_{n = N+1}^{\infty} A^n_{\alpha ,
\beta } (w_{\beta}),\cr
w_{\alpha} = w_{\beta} + \sum\limits_{n = N+1}^{\infty} B^n_{\alpha , \beta }(w_{\beta}),
\end{cases}
\end{equation}

\newprg[ROYD.psh]{Construction of a plurisubharmonic exhaustion function.}
In the previous steps we constructed a covering $\{V_{\alpha}\}$ of a neighborhood $V$
of $\phi(\bar D)$ by coordinate charts with coordinate mappings $h_{\alpha}:V_{\alpha}\to
V_{\alpha}'$ such that $V_{\alpha}'=U_{\alpha}\times V_{\alpha}^{''}$, where $U_{\alpha}$ is
a covering of $\bar D$ with coordinates $z_{\alpha}$ and $V_{\alpha}^{''}$ are open sets
in $\el^2$ with coordinates $w_{\alpha}$. And, finally, the coordinate changes have the form
\eqqref(chngN+1). We take $N=3$ in the sequel, \ie coordinates match to order three. Let 
$\rho_{\alpha}$ be nonnegative
$\calc^{\infty}$ functions with
support contained in $V_{\alpha}$ such that $\sum_{\alpha}\rho_{\alpha} = \adyn$ in a
neighborhood of $\phi(\bar D)$. Define the following (vector valued) functions on the
manifold $\calx$
\begin{equation}
\eqqno(uv1)
\begin{cases}
u = \sum\limits_{\alpha} \rho_{\alpha} z_{\alpha},\cr
v = \sum\limits_{\alpha} \rho_{\alpha} w_{\alpha}.
\end{cases}
\end{equation}
In $V_{\alpha}$ we have
\begin{equation}
\eqqno(uv2)
\begin{cases}
u - z_{\alpha} = \sum\limits_{\beta} \rho_{\beta} (z_{\beta} - z_{\alpha}) = O(w_{\alpha}^4),\cr
v  - w_{\alpha}= \sum\limits_{\beta} \rho_{\beta} (w_{\beta} - w_{\alpha}) = O(w_{\alpha}^4).
\end{cases}
\end{equation}
Since $D$ is a strongly pseudoconvex domain, it exists a strictly plurisubharmonic function $\theta$ defined in a neighborhood of $\bar{D}$ such that 
$D=\{z \in \mathbb{C}^q \ | \ \theta(z)<1 \}$. Let $\mu>0$ be the smallest eigenvalue of the Hessian matrix of $\theta$.             
Since $u$ and $v$ are differentiable coordinates in a neighborhood $\mathcal{O}$ of $\phi(\bar{D})$ we
can shrink $\mathcal{O}$ in such a way that the image of $\mathcal{O}$  under $(u,v)$ is $W\times 
B^\infty (\delta)$ with $W$ a neighborhood of $\bar{D}$. Set $\xi^2=\|v\|^2=\sum |v_k|^2.$ Then in $V_\alpha\cap\mathcal{O}$ 
we have due to \eqqref(uv2)
\[
|u_k-z_{\alpha,k}|=O(\xi^4), \qquad |v_k|^2-|w_{\alpha,k}|^2=O(\xi^5).
\]
Indeed, the first one is clear and the second one comes from 
\begin{equation}
|v_k|^2=|w_{\alpha,k}+O(w_{\alpha,k}^4)|^2=|w_{\alpha,k}|^2 +|O(w_{\alpha,k}^5)|+|O(w_{\alpha,k}^8)|=|w_{\alpha,k}|^2+O(\xi^5).
\end{equation}
For a given $0<\lambda<1$ consider the following function
\begin{equation}
\psi_\lambda (u,v) = \theta(u)+\lambda^{-2}||v||^2 = \theta (u)+\lambda^{-2}\sum\limits_k|v_k|^2.
\end{equation}
Then on $V_\alpha\cap \mathcal{O},$ we have
$$\psi_\lambda (u,v)-\psi_\lambda (z_\alpha, w_\alpha)=O(|u_k-z_{\alpha,k}|)+\lambda^{-2}
O(\sum_{k}{|v_k|^2-|w_{\alpha,k}|^2}).$$
Therefore the Hessian of the difference with respect to the coordinates $(z_\alpha, w_\alpha)$ consists of terms
$O(\xi^2)+\lambda^{-2} O(\xi^3).$ The Hessian $\mathcal{H}_\lambda$ of  $\psi_\lambda (z_\alpha, w_\alpha)$ is an 
infinite matrix of the form 
$\begin{pmatrix}
\mathcal{H}_\theta & 0 \\
0 & \lambda^{-2}I_{\infty}
\end{pmatrix}$
where $\mathcal{H}_\theta$ is the Hessian of $\theta$ and $I_{\infty}$ is the identity matrix.
The Hessian $\mathcal{H}_\lambda$ defines a positive Hermitian form larger than
$H=\begin{pmatrix}
\mu I & 0 \\
0 & \lambda^{-2}I_{\infty}
\end{pmatrix} 
$. 
Let $\mathcal{O}_\lambda$ be a subset of $\mathcal{O},$ where $\psi_\lambda (u,v)<1.$ Then 
$\bar{\mathcal{O}}_\lambda\subset \mathcal{O}$ if $\lambda<\delta.$ 
In fact, for $m \in \mathcal{O}$ inequalities $\lambda<\delta$ and $\psi_\lambda(u,v)<1$ gives $\lambda^{-2}||v(m)||^2<1$ 
and that implies $||v(m)||<\lambda<\delta$. The same goes for $\theta(u(m))<1$ that gives $u(m) \in D$. So, $\bar{\mathcal{O}}_\lambda\subset \mathcal{O}$.

In $\bar{\mathcal{O}}_\lambda\cap \bar{U}_\alpha,$ the 
Hessian of $\psi_\lambda (u,v)$ with respect to $(z_\alpha, w_\alpha)$ is a matrix which is greater than $H+O(\xi^2)
+\lambda^{-2}O(\xi^3).$ But in $\bar{\mathcal{O}}_\lambda$ we have $\xi<\lambda$ and so the Hessian is greater than
$H+O(\lambda).$ Thus there is a $\lambda_\alpha$ such that $\psi_\lambda (u,v)<1$ is strictly plurisubharmonic in
$\bar{\mathcal{O}_\lambda}\cap U_\alpha$ for $\lambda<\lambda_\alpha.$ Choose $\lambda=\min_\alpha \lambda_\alpha.$
Then $\psi_\lambda (u,v)<1$ is strictly plurisubharmonic on $\bar{\mathcal{O}}_\lambda.$ So $\eta(u,v)=
(1-\psi_\lambda (u,v))^{-1}$ is a plurisubharmonic $\calc^{\infty}$ function on $\bar{\mathcal{O}}_\lambda.$

\smallskip\qed

\newsect[HART]{Hilbert-Hartogs manifolds}

\newprg[HART.hart]{Hartogs manifolds}
 In introduction
we gave the notion of a $q$-Hartogs manifold. $1$-Hartogs manifolds will be called simply 
Hartogs. A closed submanifold of a $q$-Hartogs
manifold is $q$-Hartogs. If, for example, $P(z) = \sum_{k=1}^{\infty}z_k^2$ then the
manifold $\calx = \{z\in \el^2: P(z) = 1\}$ is Hartogs and analogously any $\calx$, which is defined
as a non singular zero set of finitely many continuous polynomials in $\el^2$ is Hartogs.

\smallskip The following statement is immediate via the Docquier-Grauert characterization
of Stein domains over Stein manifolds, see \cite{DG}.

\begin{prop}
\label{doq-gra}
If $\calx$ is a Hartogs manifold then for every domain $D$
over a Stein manifold every holomorphic mapping $f:D\to \calx $
extends to a holomorphic mapping $\hat f:\hat D\to \calx$ of the
envelope of holomorphy $\hat D$ of $D$ to $\calx$.
\end{prop}

As we told in the introduction Hartogs manifolds have better
extension properties than it is postulated in their definition. Let
us make a first step in proving this. For positive integers $q$, $n$
and real $0< r <1$ we call a Hartogs figure of bidimension $(q,n)$
or, a $q$-concave Hartogs figure  in $\cc^{q+n}$ the set

\begin{equation}
\eqqno(hart-qn)
H_q^n(r):= \left(\Delta^q\times \Delta^n(r)\right)\cup \left(A^q_{1-r ,1}\times \Delta^n\right).
\end{equation}
The envelope of holomorphy of $H_q^n(r)$ is $\Delta^{q+n}$.
A $q$-concave Hartogs figure $H_q^1$ from the introduction is a Hartogs figure of bidimension $(q,1)$.
If holomorphic maps with values in Hilbert manifold $\calx$
holomorphically extend from $H_q^n(r)$ to $B^q\times
B^n$ we shall say that $\calx$ possesses a holomorphic
extension property in bidimension $(q,n)$. Let us prove that
holomorphic extendability in bidimension $(q,n)$, \ie the property
being $q$-Hartogs, implies the
holomorphic extendability in all bidimensions $(p,m)$ with $p\ge q$
and $m\ge n$. The proof follows the lines of Lemmas 2.2.1 and 2.2.2
in \cite{Iv2}, \ie the finite dimensional case and therefore we
shall give only the detailed proof of the first step, pointing out
at what place one needs Theorem \ref{royden-w} for Hilbert manifolds.
All other steps are just repetitions of the first one and will be
only sketched.

\begin{rema} \rm
\label{ball-disk} When working with Hartogs figures one can replace
$B^q$ (resp. $B^n$) in their definitions by the polydisks of the
corresponding dimension. This is not the case when one works
with $B^{\infty}$ - the unit ball in $\el^2$. 
\end{rema}

\begin{thm}
\label{qn-pm} If a complex Hilbert manifold $\calx$ possesses a
holomorphic extension property in bidimension $(q,n)$ for some $q, n
\ge 1$ then $\calx$ possesses this property in every bidimension
$(p,m)$ with $p\ge q, m\ge n$.
\end{thm}
\proof We shall prove this by induction. Let us increase $q$ first. Remark that
\[
 H_{q+1}^n(r) \supset \bigcup_{z_{q+1}\in \Delta }  H_q^n(r)\times \{z_{q+1}\},
\]
and therefore a holomorphic mapping $f: H_{q+1}^n(r)\to \calx$
extend along these slices to a map $\tilde f : \Delta^{q}\times
\Delta \times \Delta^{n}\to \calx$, we respect here the order of
variables: $z\in \cc^q$, $z_{q+1}\in \cc$, $w\in\cc^n$ and use the
notation $Z = (z,z_{q+1},w)$ for these coordinates. 

\smallskip\noindent{\sl Continuity of $\tilde f$.} First we need to prove
that this extension $\tilde f$ is continuous. Let a sequence $Z_k
= (z_1^k,...,z_q^k,z_{q+1}^k,w_1^k, ...,w_n^k)$ converge to $Z_0 =
(z_1^0,...,z_q^0,z_{q+1}^0,w_1^0, ...,w_n^0)\in \Delta^{q+1+n}$.
Take $0<R<1$ such that $\norm{z^0}, |z_{q+1}^0|, \norm{w^0}<R$ and
the same is true for $\norm{z^k}, |z^k_{q+1}|, \norm{w^k}$ with $k$
big enough. Consider the following imbedded $(q+n)$-disk $\Phi^0$ in
the Hilbert manifold $\caly \deff \cc^{q+1+n}\times \calx$:
\[
\phi_0 (z,w) = \Big\{\left(z, z_{q+1}^0, w, \tilde f(z,z_{q+1}^0,w)\right): z\in \bar\Delta^q(R),
w\in\bar\Delta^n(R)\Big\}
\]
\ie $\Phi^0$ is the graph of the restriction of $\tilde f$ to the $(q+n)$-disk $\bar\Delta^q(R)\times
\{z^0_{q+1}\}\times \bar\Delta^{n}(R) \subset \Delta^{q+1+n}$.

\smallskip It will be convenient for the future references to formulate the next
step of the proof in the form of a lemma.

\begin{lem}
\label{cont}
Let a $\phi_n:\bar D\to \calx$ be a sequence of analytic disks in a Hilbert manifold $\calx$ and
let $\Phi_n$ be their graphs, here $D\comp \cc^q$. Suppose that there exists an analytic disk 
$\phi_0:\bar D\to \calx$ with the graph $\Phi_0$ such that for any neighborhood $\calv\supset \Phi_0$ 
one has $\Phi_n\subset \calv$ for $n>>1$. Then $\phi_n$ converge uniformly on $\bar D$ to $\phi_0$.
\end{lem}
\proof  Since $\phi_0$ is uniformly continuous on $\bar D$ for a given $\eps >0$ we can find $\delta >0$
such that if $||u - v|| <\delta$ then $d(\phi_0(u),\phi_0(v))< \eps$. Here $d$ is some Riemannian
metric on $\calx$. In addition we can assume that $\delta <\eps$. 
Take $N$ such that for $\forall n\ge N$ one has  $\Phi_n\subset \Phi_0^{\delta}$,
where $\Phi_0^{\delta}$ is a $\delta$-neighborhood of $\Phi_0$ with respect to the product
metric on $\calz \deff \cc^q\times \calx$. Fix $u\in \bar D$. For every $n\ge N$ there exists $v_n\in \bar D$ 
such that the distance between $(u, \phi_n(u))$ and $(v_n,\phi_0(v_n))$ is less than $\delta$. But then the 
distance between $u$ and $v_n$ is less than $\delta$, \ie $||u-v_n||<\delta$. Therefore $d(\phi_0(u),\phi_0(v_n))<\eps$.
This proves that the distance between $(u, \phi_n(u))$ and $(u,\phi_0(u))$ is not more than $\eps + \delta
<2\eps$.

\smallskip\qed

\smallskip To apply this lemma in our setting remark that
by Theorem \ref{royden-w} there exists a neighborhood $\calv$ of $\Phi^0$ in
$\caly\deff D^r\times\calx$  and a plurisubharmonic function $\rho:\bar \calv\to \rr^+$ such that
$\calv_r\deff \{\rho <r\}$ decreases to $\Phi_0$. Remark furthermore that for a given $r>0$ for $z_{q+1}$ close
enough to $z_{q+1}^0$
the graph of $\tilde f$ over  of $H^n_q\times\{z_{q+1}\}$ belongs to $\calv_r$. Therefore its graph
over the whole polydisk $\bar\Delta^{n+q}(R)\times \{z_{q+1}\}\times \bar \Delta^n(R)$ belongs to
$\calv_r$ by maximum principle. Therefore the question of
continuity is reduced to Lemma \ref{cont} just proved.

\smallskip\noindent{\sl Holomorphicity of $\tilde f$.} Holomorphicity becomes now a local question and therefore 
is obvious. Indeed, let $\{ \Omega_\alpha, \phi_\alpha \}$ an atlas of $\mathcal{X}$. Cover $\Delta^q(R) \times 
\Delta(R) \times \Delta^n(R)$ by polydiscs of polyradius $\epsilon$. As it is compact, there is a finite number
of polydiscs and let note $x^j$ the j-th centers. By choosing $\epsilon$ small enough we can suppose that for 
all $x^j$ there exists $\alpha_j$ such that 
$f(\Delta(x^j, \epsilon)) \subset \Omega_{\alpha_j}$.
Let us define the mapping $ \psi_{x^j}$ by
\begin{equation}
\psi_{x^j}(z,z_{q+1},w)= \phi_{\alpha_j}^{-1}( \frac{1}{2 \pi i}\int_{\partial \Delta(x_{q+1}^j,\epsilon)}{\frac{\phi_{\alpha_j} 
\circ f(z,\xi,w)}{\xi - z_{q+1}}d\xi} \ )
\end{equation}

$ \psi_{x^j}$ are by separable analyticity $(z,w)$-holomorphic and are $z_{q+1}$ holomorphic by Cauchy formula. 
The Hartogs theorem gives that $ \psi_{x^j}$ are (globally) holomorphic.Also, $ \psi_{x^j}$ coincide with 
$ \psi_{x^l}$ ($j \neq l$) by unicity property of holomorphic function. Then, we can define a global function 
$ \Psi$ on $\Delta^q(R) \times \Delta(R) \times \Delta^n(R)$. As 
$\Psi$ coincide with $\tilde{f}$ on $H_{q+1}^n(r)$, one has $\Psi = \tilde{f}$. So $\tilde{f}$ is holomorphic.
$f$ holomorphically extended to $\Delta^q\times\Delta
\times\Delta^n$ as required.

\smallskip Increase of  $n$ follows the same lines and can be fulfilled also in two steps. Set
\[
E_q^{n+1}(r) = H_q^n(r)\times \Delta (r),
\]
and  remark that $E_q^{n+1}(r)\subset H_q^{n+1}(r)$. Extend $f$ from
$E_q^{n+1}(r)$ to $\Delta^{q+n}\times \Delta (r)$ exactly as above. Then extend $f$
from $\Delta^{q+n}\times \Delta (r)$ to $\Delta^{q+n+1}$ again in the same way.

\smallskip\qed

\begin{rema} \rm
\label{non-mer} In \cite{Iv2} it was shown with an example that
meromorphic extendability in bidimension $(1,1)$ doesn't imply that
in bidimensions neither $(1,2)$ (which is not surprising) but also
not in bidimension $(2,1)$, even if the image manifold is finite
dimensional.
\end{rema}

\newprg[HART.levi]{Levi extension property}

Now we shall prove Theorem \ref{hart-levi} from the Introduction.

\medskip\noindent {\sl $q$-Hartogs $\Rightarrow$ $q$-Levi extension.}
Let $H_q^1(r)=(\Delta^q \times \Delta(r)) \cup (A^q_{1-r,1} \times \Delta)$ the $q$-Hartogs figure. 
Let $\calx$ be $q$-Hartogs and consider a domain $D$ as in the definition. Let $U$ an open set in $\mathbb{C}^{n}$ such that 
$U \cap D = \{z \in U \ | \ u(z)<0 \}$ with $\nabla u \neq 0$.  Then, there exist $v_j \in T_0\partial D \subset  \mathbb{C}^n$ such that $\mathcal{L}_{u,p}(v_j)<0$ for $j=1,...,q$. After a complex linear change of coordinate, we can suppose that 
$p=0$ and 

\[
v_1=e_1, \ v_2=e_2, \ ... \ , \ v_q=e_q, \nabla u_p = e_{q+1}.
\]

We can define the following figure:

\[
\phi : H_q^1(r) \times \Delta^{n-q-1} \to \phi(H_q^1(r)\times \Delta^{n-q-1}) 
\]
\begin{equation}
z \mapsto \eta z_1 e_1 + ...+ \eta z_q e_q + (\eta \epsilon z_{q+1} -r')e_{q+1}
+\delta z_{q+2}e_{q+2}+...+\delta z_n e_n.
\end{equation}

We can choose $\epsilon$, $r'$, $\delta$ and $\eta$ such that $\phi(H_q^1(r)\times \Delta^{n-q-1}) 
\subset U \cap D$ and $\eta \epsilon >r'$. The latter insures that the image of the polydisk $\phi (\Delta^n)$ 
contains the origin.

\begin{figure}[h]
\centering
\includegraphics[width=2.0in]{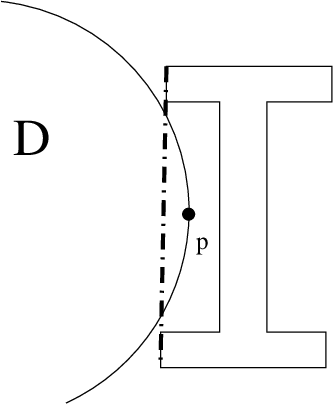}
\caption{$q$-Hartogs figure in a neighborhood of p}
\label{pic2}
\end{figure}

Now, for any $f:D \rightarrow \calx$, the map $f|_{\phi(H_q^1(r) \times \Delta^{n-q-1})}$ extends to $ \phi(\Delta^{q+1} 
\times \Delta^{n-q-1})$ because $\calx$ is $q$-Hartogs. So it extends to a neighborhood of $p=0$. Therefore $\calx$ 
satifies the $q$-Levi extension condititon.

\smallskip\noindent{\sl $q$-Levi extension $\Rightarrow$ $q$-Hartogs.}
Suppose $\calx$ satisfies the $q$-Levi extension condition and let $f:H_q^1(r) \rightarrow \calx$ be a holomorphic map. 
Consider, for $2\le \alpha < \infty$ and $0<\eps <1$  the following function

\begin{equation}
\eqqno(al-eps)
u_{\alpha , \eps}(z)=|z_{q+1}| -(1-\eps)||z'||^{\alpha} - \eps ,
\end{equation}
where $z'=(z_1,...,z_q)$ and let $D_\alpha = \{ z \in B^q \times \Delta \ | \ u_{\alpha , \eps}(z)<0 \}$.
One can see that:

\begin{itemize}

\item For a fixed $0<\eps_0<1$ here exist $\alpha_0 \ge 2$ such that $\overline{D}_{\alpha_0 , \eps_0} \subset  H_1^q(r)$.

\item $D_{\alpha , \eps}$ satisfies the $q$-Levi extension condition for all $0<\eps<1$ and $2\le \alpha <\infty $.

\end{itemize}

Let $\Gamma_{\epsilon_0} = \{\alpha \geq 2 \ | \ f $ extends to $D_{\alpha,\epsilon_0}$\}.
Then $\Gamma_{\epsilon_0}$ is nonempty (because $\alpha_0 \in \Gamma_{\epsilon_0}$) and closed.
As $X$ satisfies the levi extension condition, $\Gamma_{\epsilon_0}$ is open. 
$f$ extends on $D_{2,\epsilon_0}$. 
Now, let $\Gamma = \{\epsilon <1 \ | \ f $ extends to $D_{2,\epsilon}$\}. For the same reasons above, $\Gamma$
is nonempty and closed. It is also open because $\mathcal{X}$ satisfies the $q$-Levi extension condition. Hence 
$f$ extends holomorphically on $D_{2,1}$: 
\begin{equation}
D_{2,1}= \{ z \in \Delta^{q+1} \ | \ |z_{q+1}|-1<0 \}= \Delta^{q+1}
\end{equation}
$f$ extends from $H_q^1(r)$ to $\Delta^{q+1}$. Then $\mathcal{X}$ is $q$-Hartogs.

\smallskip\qed

\smallskip
After the preparations already made the following statement can be proved along the same lines as 
Proposition 4 in \cite{Iv1}.

\begin{prop}
\label{fiber}
(a) If $\calx$ is a Hilbert manifold and $\caly$ is some unramified
cover of $\calx$ then $\calx$ and $\caly$ are $q$-Hartogs or not simultaneously.

(b) If the fiber $\calf$ and the base $\calb$ of a complex Hilbert fiber bundle
$(\cale,\calf,\pi ,\calb)$ are $q$-Hartogs then the total space $\cale$ is also
$q$-Hartogs.
\end{prop}
\proof (a) Suppose that $\mathcal{X}$ is $q$-Hartogs. We shall show that $\mathcal{Y}$ satisfies the Levi extension condition:
For $D$ a domain in $\mathbb{C}^{q+1}$ with $\mathcal{C}^2$ boundary  such that $\partial D$ is $q$-Levi pseudoconcave at $p$.
Let $f:D \rightarrow \mathcal{Y}$ be holomorphic. 
As $\mathcal{X}$ is $q$-Hartogs, we can extend $\pi \circ f$ to an open neighborhood $V$ of $p$ with $\pi$ the projection map.
Take a neighborhood $U \subset V$ of $p$ such that there is a trivialisation: 
$$h=(\pi,b): \pi^{-1}(U) \tilde{\rightarrow} U \times F$$
As $F$ is a discrete set, so we can extend $f$ on a neighborhood of $h^{-1}(U \times \{b\circ f(p)\}) \subset \mathcal{Y}$ 
Hence, $f$ satisfies the $q$-Levi extension condition and is $q$-Hartogs.

Reciprocally, suppose that $\mathcal{Y}$ is $q$-Hartogs. For $D$ a domain in $\mathbb{C}^{q+1}$ with $\mathcal{C}^2$ boundary 
such that $\partial D$ is $q$-Levi pseudoconcave at $p$.
Let $f:D \rightarrow \mathcal{X}$ be holomorphic.
For $p \in D$, it exsits $\mathfrak{p} \in \mathcal{Y}$ such that $\pi(\mathfrak{p})=f(p)$. We take a neighborhood $V$ of 
$\mathfrak{p}$ that is biholomorphic to $\pi(V)$ because $\pi$ is locally a biholomorphism. Then, $(\pi|_V)^{-1}\circ f$ 
extends to a neighborhood of $p$. Hence, we compose the last extension by $\pi$ and it gives an extension of $f$ to a 
neighborhood of $p$.
$\mathcal{X}$ satisfies the $q$-Levi extension condition so is $q$-Hartogs.

\smallskip \noindent (b) We shall show that $\mathcal{E}$ satisfies the Levi extension condition.
For $D$ a domain in $\mathbb{C}^{q+1}$ with $\mathcal{C}^2$ boundary  such that $\partial D$ 
is $q$-Levi pseudoconcave at $p$.
Let $f:D \rightarrow \mathcal{E}$ be holomorphic. 
As $\mathcal{B}$ is $q$-Hartogs, we can extend $\pi \circ f$ to an open neighborhood $V$ of $p$.
Take a neighborhood $U \subset V$ of $p$ such that there is a trivialisation: 
$$h=(\pi,p): \pi^{-1}(U) \tilde{\rightarrow} U \times \mathcal{F}$$
Then, $p \circ f$ extends to a neighborhood of $p$ as $\mathcal{F}$ is $q$-Hartogs. So $\mathcal{E}$ 
satisfies the $q$-Levi extension condition and is $q$-Hartogs.

\smallskip\qed

\begin{rema} \rm From this proposition one can derive several examples of Hilbert-Hartogs manifolds.
If, for example, $\Lambda = \span_{\zz}\{e_1, ie_1, e_2, ie_2,...\}$ is the integer lattice in
$\el^2$, then $\ttt^{\infty} \deff \el^2/\Lambda$ is Hartogs. Hilbert fiber bundles over Riemann surfaces
different from $\pp^1$ are Hilbert-Hartogs.
\end{rema}

\ifx\undefined\bysame
\newcommand{\bysame}{\leavevmode\hbox to3em{\hrulefill}\,}
\fi

\def\entry#1#2#3#4\par{\bibitem[#1]{#1}
{\textsc{#2 }}{\sl{#3} }#4\par\vskip2pt}

\end{document}